\documentclass[11pt]{article}
\usepackage[letterpaper,margin=1in]{geometry}
\usepackage[T1]{fontenc}
\usepackage[utf8]{inputenc}
\usepackage{lmodern}
\usepackage{microtype}
\usepackage{amsmath,amssymb,amsthm,mathtools}
\usepackage{booktabs}
\usepackage{enumitem}
\usepackage[hidelinks]{hyperref}
\usepackage{xurl}
\usepackage{array}
\usepackage{longtable}

\newtheorem{theorem}{Theorem}[section]
\newtheorem{lemma}[theorem]{Lemma}
\newtheorem{proposition}[theorem]{Proposition}
\newtheorem{corollary}[theorem]{Corollary}
\theoremstyle{remark}
\newtheorem*{remark}{Remark}

\DeclareMathOperator{\Logop}{Log}
\DeclareMathOperator{\arccsc}{arccsc}
\DeclareMathOperator{\arcsec}{arcsec}
\DeclareMathOperator{\arsinh}{arsinh}
\DeclareMathOperator{\sech}{sech}
\DeclareMathOperator{\sgn}{sgn}
\DeclareMathOperator{\gd}{gd}

\title{A Unified Substitution Method for Integration (DRAFT)}
\author{Emmanuel Antonio José García\\CIDIC-UTE\\Dominican Republic\\emmanuelgeogarcia@gmail.com}
\date{August 22, 2026}

\begin{document}
\maketitle

\begin{abstract}
We present a branch-consistent framework for integrals involving quadratic radicals by expressing exponentials of principal inverse trigonometric functions in algebraic form. Two identities for
\[
e^{\pm i\arccos(y)}
\qquad\text{and}\qquad
e^{\pm i\arcsec(y)}
\]
on principal branches yield five explicit substitution templates that map common radicals and half-angle composites to rational functions of a single parameter. The resulting differentials are independent of the sign choice once the branches are fixed, reducing domain bookkeeping across circular and hyperbolic regimes. For several natural families, cancellation between the parametrized radical and the Jacobian makes the transformed integrand a finite Laurent polynomial, so it can be integrated term by term. We recover Euler's first and second substitutions from these transforms
after the necessary scaling, sign choices, and component-dependent
reciprocal reparametrizations and provide worked examples; in particular,
the classical Weierstrass substitution is obtained as a local direct
corollary of Transform~5. The framework also yields Lobachevsky's formula
\[
\tan\!\left(\frac{\Pi(u)}{2}\right)=e^{-u}
\]
as the positive hyperbolic specialization of one core identity, with an equivalent derivation through the reciprocal coordinate \(y=\sech u\). Binomial sum-and-difference identities streamline the back-substitution of paired Laurent powers such as \(p^n\pm p^{-n}\). Under the reported computational protocol, all tested instances for Transforms~2--5 passed the reported derivative/numerical checks, as did 94 of 100 instances for Transform~1. Selected benchmark families also exhibited reduced expression swell or improved completion coverage relative to the native Wolfram Language function \texttt{Integrate} under the reported protocol, although these results concern the targeted test families and do not establish general superiority.
\end{abstract}

\section{Introduction}

Classical techniques for integrating expressions involving quadratic
radicals, including
\[
\sqrt{x^2-a^2},\qquad \sqrt{a^2-x^2},\qquad \sqrt{x^2+a^2},
\]
use trigonometric, hyperbolic, or Euler substitutions
\cite{Hardy,Kishan,Piskunov}. On real domains, the resulting
parametrizations must be interpreted on appropriate components with
the corresponding radical signs fixed.The purpose of this article is organizational: the Unified Substitution Method (USM) expresses several of these substitutions, together with related inverse-trigonometric half-angle formulas, in a common parameter language while retaining the branch data required for correct back-substitution.

The core identities were originally discovered in the course of an attempt to generalize the Newton--Mollweide formula \cite{GarciaMollweide}. That investigation led to a family of trigonometric--hyperbolic representations for the roots of quadratic equations \cite{GarciaRoots}, and the identities used here emerged by equating corresponding root formulas. In the present paper, however, they are rederived independently and more directly from Euler's formula together with the adopted principal-branch definitions of the inverse trigonometric functions. They provide a compact organization of a useful family of substitutions and retain a classical geometric specialization. Section~\ref{sec:lobachevsky} shows that Theorem~\ref{thm:cos-csc}, on the positive hyperbolic component $y=\cosh u$, yields Lobachevsky's formula $\tan(\Pi(u)/2)=e^{-u}$; Theorem~\ref{thm:sec-sin} yields the same formula through the reciprocal coordinate $y=\sech u$.

In some cases, the practical advantage is stronger than rationalization alone. When a parametrized radical and its Jacobian share factors, cancellation occurs before integration, and the transformed integrand may become a finite Laurent polynomial
\[
\sum_{k=m}^{n}c_kp^k,
\qquad m,n\in\mathbb Z,\quad m\le n,
\]
where $p$ is the relevant USM parameter. This occurs systematically for several natural families. With $X=x+b$, examples include
\[
P(x)\sqrt{\frac{X-a}{X+a}}\,dx,
\qquad
\frac{dx}{X^n\sqrt{a^2-X^2}},
\qquad
\frac{1}{X^n}\sqrt{\frac{a+X}{a-X}}\,dx,
\]
subject to the index restrictions in Proposition~\ref{prop:laurent}. Transform~3 provides further systematic Laurent families. In particular, for every polynomial $P$, both
\[
P(x)\sqrt{X^2+a^2}\,dx,
\qquad
\frac{P(x)}{\sqrt{X^2+a^2}}\,dx
\]
have Laurent-polynomial transformed forms, as do reciprocal integer powers of $X+\sqrt{X^2+a^2}$ and $\sqrt{X^2+a^2}-X$. When $b=0$, powers of $X$ are simply powers of $x$. Such integrals can be evaluated term by term. Examples~1, 3, 4, and 6 illustrate this outcome.

Throughout, transformed antiderivative formulas are asserted on open connected intervals contained in the interior of the displayed real domain and avoiding every zero of a displayed denominator. Some algebraic parametrization formulas also have well-defined endpoint values. Such values are used in definite integrals or endpoint extensions only when the complete expression has a verified one-sided limit. Constants of integration may depend on the connected component.

Appendix~\ref{app:continuation} records complementary half-angle continuations and principal-root sign tables, and Appendix~\ref{app:binomial} collects binomial sum-and-difference formulas for reconstructing paired Laurent powers. The public post \emph{Benchmarking the USM Against Mathematica Integrate, Part 3} \cite{GarciaBenchmark} links a computational report containing the displayed Wolfram Language benchmark material; Appendix~\ref{app:benchmark} summarizes this author-reported evidence and its limitations.

\section{Principal branches and boundary-value conventions}
\label{sec:branches}

Let $\Logop$ denote the analytic principal logarithm on
\[
\mathbb C\setminus(-\infty,0],
\qquad -\pi<\arg z<\pi.
\]
On this slit plane, define
\[
\sqrt z:=\exp\!\left(\frac{1}{2}\Logop z\right),
\]
and extend the function continuously to the origin by $\sqrt0:=0$. When a displayed square-root radicand lies on the negative real axis, we use the upper-bank boundary value
\[
\sqrt{-\rho+i0}:=\lim_{\varepsilon\downarrow0}\sqrt{-\rho+i\varepsilon}=i\sqrt\rho,
\qquad \rho>0.
\]
For brevity, the notation $\sqrt{-\rho}=i\sqrt\rho$ below always denotes this boundary value. It is not an assertion that the analytic square root is holomorphic on its branch cut. The phrase ``upper bank'' refers to the square-root radicand; after a change of variable or reciprocal composition, it need not correspond to a fixed bank in the original variable. These conventions agree with standard treatments of principal branches \cite{Ahlfors,BrownChurchill}.

Define
\[
\arcsin z:=-i\Logop\!\left(iz+\sqrt{1-z^2}\right),
\qquad
\arccos z:=\frac{\pi}{2}-\arcsin z.
\]
These are analytic principal branches on
\[
\Omega_{\sin/\cos}:=\mathbb C\setminus\bigl(( -\infty,-1]\cup[1,\infty)\bigr).
\]
For $z\ne0$, define
\[
\arcsec z:=\arccos(1/z),
\qquad
\arccsc z:=\arcsin(1/z).
\]
Because the reciprocal map pulls the cuts of $\arccos$ and $\arcsin$ back to $[-1,1]$, these functions are analytic principal branches on
\[
\Omega_{\sec/\csc}:=\mathbb C\setminus[-1,1].
\]
Thus the two cut planes are different, in agreement with standard principal-branch conventions for inverse trigonometric functions \cite{DLMF}.

At $z=\pm1$, all four functions are assigned their continuous endpoint values. These values agree with the displayed formulas when $\sqrt0=0$, although $\pm1$ remain branch points and do not belong to the corresponding analytic cut-plane domains. On the real axis, the ordinary principal values satisfy
\[
-1\le y\le1:
\quad
\arccos y\in[0,\pi],
\quad
\arcsin y\in[-\pi/2,\pi/2],
\]
and
\[
|y|\ge1:
\quad
\arcsec y\in[0,\pi]\setminus\{\pi/2\},
\quad
\arccsc y\in[-\pi/2,0)\cup(0,\pi/2].
\]
The endpoint values at $y=\pm1$ are understood by continuity. The functions $\arcsec0$ and $\arccsc0$ are undefined.

For real $|y|>1$, the symbols $\arcsin y$ and $\arccos y$ below denote the boundary values obtained from the displayed logarithmic formulas and the preceding upper-bank square-root convention. For real $0<|y|<1$, the symbols $\arcsec y$ and $\arccsc y$ denote the corresponding boundary values inherited through
\[
\arcsec y=\arccos(1/y),
\qquad
\arccsc y=\arcsin(1/y).
\]
Accordingly,
\[
\sqrt{y^2-1}=i\sqrt{1-y^2},\qquad |y|<1,
\]
and
\[
\sqrt{1-y^2}=i\sqrt{y^2-1},\qquad |y|>1.
\]
See \cite[\S4.4]{AbramowitzStegun} for the standard real inverse-function ranges.

\section{Core identities}

Throughout this section, all inverse trigonometric functions and square
roots are interpreted according to the principal-branch and boundary-value
conventions of Section~\ref{sec:branches}.

\begin{theorem}[Cosine--cosecant identity]
\label{thm:cos-csc}
Let \(y\in\mathbb R\) and set
\[
\alpha=\arccos y.
\]
When \(y\ne0\), also set
\[
\beta=\arccsc y,
\qquad
\gamma=\arcsec y.
\]
Then:
\begin{enumerate}[label=\arabic*.]
\item If \(|y|\le1\), then
\[
e^{-i\alpha}
=y-\sqrt{y^2-1}
=y-i\sqrt{1-y^2}.
\]

\item If \(|y|\ge1\), then
\[
\tan\!\left(\frac{\beta}{2}\right)
=
\begin{cases}
e^{i\alpha}
=y-\sqrt{y^2-1},
& y\ge1,\\[2mm]
e^{-i\alpha}
=y+\sqrt{y^2-1},
& y\le-1.
\end{cases}
\]
Whenever \(\tan(\gamma/2)\) is defined,
\[
\tan\!\left(\frac{\beta}{2}\right)
=
\frac{1-\tan(\gamma/2)}
     {1+\tan(\gamma/2)}.
\]
At \(y=-1\), the last equality holds only as a one-sided limit because
\(\tan(\gamma/2)=\tan(\pi/2)\) is undefined.
\end{enumerate}
\end{theorem}

\begin{proof}
For \(|y|\le1\), one has \(\alpha\in[0,\pi]\), and hence
\[
e^{-i\alpha}
=\cos\alpha-i\sin\alpha
=y-i\sqrt{1-y^2}
=y-\sqrt{y^2-1}.
\]

For \(|y|\ge1\), the half-angle identity
\[
\tan\!\left(\frac{\beta}{2}\right)
=\frac{1-\cos\beta}{\sin\beta}
\]
together with
\[
\sin\beta=\frac{1}{y},
\qquad
\cos\beta
=\sqrt{1-\frac{1}{y^2}}
=\frac{\sqrt{y^2-1}}{|y|}
\]
gives
\[
\tan\!\left(\frac{\beta}{2}\right)
=
\begin{cases}
y-\sqrt{y^2-1},
& y\ge1,\\[1mm]
y+\sqrt{y^2-1},
& y\le-1.
\end{cases}
\]
By the standing conventions,
\[
e^{\pm i\alpha}
=y\mp\sqrt{y^2-1},
\]
so the stated sign selection follows.

Finally, since
\[
\beta=\arcsin(1/y),
\qquad
\gamma=\arccos(1/y),
\]
the defining relation \(\arccos z=\pi/2-\arcsin z\) gives
\[
\beta+\gamma=\frac{\pi}{2}.
\]
Therefore, wherever \(\tan(\gamma/2)\) is defined,
\[
\tan\!\left(\frac{\beta}{2}\right)
=\tan\!\left(\frac{\pi}{4}-\frac{\gamma}{2}\right)
=\frac{1-\tan(\gamma/2)}
     {1+\tan(\gamma/2)}.
\]
At \(y=-1\), the right-hand side is interpreted by the stated one-sided limit.
\end{proof}

\begin{theorem}[Secant--sine identity]
\label{thm:sec-sin}
Let \(y\in\mathbb R\) and set
\[
\psi=\arcsin y,
\qquad
\alpha=\arccos y.
\]
When \(y\ne0\), also set
\[
\gamma=\arcsec y.
\]
Then:
\begin{enumerate}[label=\arabic*.]
\item If \(0<y\le1\), then
\[
e^{i\gamma}
=\tan\!\left(\frac{\psi}{2}\right).
\]
If \(-1\le y<0\), then
\[
e^{-i\gamma}
=\tan\!\left(\frac{\psi}{2}\right).
\]
At \(y=0\), \(\arcsec 0\) is undefined, but
\[
\tan\!\left(\frac{1}{2}\arcsin 0\right)=0.
\]

\item If \(y\ge1\), then
\[
e^{-i\gamma}
=\tan\!\left(\frac{\psi}{2}\right),
\]
and if \(y\le-1\), then
\[
e^{i\gamma}
=\tan\!\left(\frac{\psi}{2}\right).
\]
Equivalently, for \(|y|\ge1\),
\[
\tan\!\left(\frac{\psi}{2}\right)
=\frac{1-i\sqrt{y^2-1}}{y}.
\]

\item On \(-1<y\le1\),
\[
\tan\!\left(\frac{\psi}{2}\right)
=
\frac{1-\tan(\alpha/2)}
     {1+\tan(\alpha/2)}.
\]
At \(y=-1\), the right-hand side is undefined and has limiting value
\(-1\).
\end{enumerate}
\end{theorem}

\begin{proof}
For \(0<y\le1\), write
\[
y=\sech u,
\qquad
u\ge0.
\]
Then \(\gamma=iu\), and therefore
\[
\tan\!\left(\frac{\psi}{2}\right)
=\frac{1-\sqrt{1-y^2}}{y}
=\cosh u-\sinh u
=e^{-u}
=e^{i\gamma}.
\]

For \(-1\le y<0\), write
\[
y=-\sech u.
\]
Then \(\gamma=\pi-iu\), and the same calculation gives
\[
\tan\!\left(\frac{\psi}{2}\right)
=-e^{-u}
=e^{-i\gamma}.
\]
The value at \(y=0\) follows from \(\psi=0\).

For \(|y|\ge1\),
\[
\sin\psi=y,
\qquad
\cos\psi=\sqrt{1-y^2}
=i\sqrt{y^2-1},
\]
so
\[
\tan\!\left(\frac{\psi}{2}\right)
=\frac{1-i\sqrt{y^2-1}}{y}.
\]
Because
\[
\cos\gamma=\frac{1}{y}
\]
and the real principal value of \(\gamma\) has nonnegative sine,
\[
e^{\pm i\gamma}
=\frac{1}{y}
 \pm i\frac{\sqrt{y^2-1}}{|y|},
\]
and comparison gives the stated componentwise signs.

Finally, on \(-1<y\le1\), the defining relation
\[
\alpha=\arccos y=\frac{\pi}{2}-\arcsin y=\frac{\pi}{2}-\psi
\]
gives
\[
\frac{\psi}{2}=\frac{\pi}{4}-\frac{\alpha}{2}.
\]
Hence
\[
\tan\!\left(\frac{\psi}{2}\right)
=\tan\!\left(\frac{\pi}{4}-\frac{\alpha}{2}\right)
=\frac{1-\tan(\alpha/2)}
     {1+\tan(\alpha/2)}.
\]
This calculation is valid at \(y=1\); at \(y=-1\), only the limit is
asserted.
\end{proof}

\section{Lobachevsky's formula as a geometric specialization}
\label{sec:lobachevsky}

Restrict Theorem~\ref{thm:cos-csc} to its positive hyperbolic component by setting
\[
y=\cosh u,
\qquad u\ge0,
\]
and hence
\[
\alpha=\arccos(\cosh u),
\qquad
\beta=\arccsc(\cosh u).
\]
The Gudermannian satisfies
\[
\gd(u)=\arctan(\sinh u)
\]
and therefore
\[
\cos(\gd(u))=\frac{1}{\sqrt{1+\sinh^2u}}=\sech u
\]
\cite[\S4.23(viii)]{DLMF}. It follows that
\[
\sin\!\left(\frac{\pi}{2}-\gd(u)\right)=\sech u.
\]
On the other hand,
\[
\sin\beta=\frac{1}{\cosh u}=\sech u.
\]
Because both $\beta$ and $\pi/2-\gd(u)$ belong to $(0,\pi/2]$, we obtain
\[
\beta=\frac{\pi}{2}-\gd(u).
\]
In the standard curvature $-1$ normalization, Lobachevsky's angle of parallelism $\Pi(u)$ is classically complementary to the Gudermannian:
\[
\Pi(u)+\gd(u)=\frac{\pi}{2}
\]
\cite[\S16.3]{Coxeter}. Consequently,
\[
\beta=\arccsc(\cosh u)=\Pi(u).
\]

It remains to specialize the exponential identity in Theorem~\ref{thm:cos-csc}. Under the boundary-value convention of Section~\ref{sec:branches},
\[
\sqrt{1-\cosh^2u}=i\sinh u.
\]
Therefore,
\begin{align*}
\arcsin(\cosh u)
&=-i\Logop\!\left(i\cosh u+\sqrt{1-\cosh^2u}\right)\\
&=-i\Logop\!\left(i(\cosh u+\sinh u)\right)\\
&=-i\Logop(ie^u)=\frac{\pi}{2}-iu,
\end{align*}
and hence
\[
\alpha=\arccos(\cosh u)=iu.
\]
The positive exterior component of Theorem~\ref{thm:cos-csc} now gives
\[
\tan\!\left(\frac{\Pi(u)}{2}\right)
=\tan\!\left(\frac{\beta}{2}\right)
=e^{i\alpha}=e^{-u}.
\]
Thus Lobachevsky's classical formula
\[
\tan\!\left(\frac{\Pi(u)}{2}\right)=e^{-u}
\]
is precisely the positive hyperbolic geometric specialization of Theorem~\ref{thm:cos-csc} \cite[\S36]{Lobachevsky}. The same conclusion follows explicitly from Theorem~\ref{thm:sec-sin}: with $y=\sech u$, one has $\psi=\arcsin(\sech u)=\Pi(u)$ and $\gamma=\arcsec(\sech u)=\arccos(\cosh u)=iu$, so $\tan(\Pi(u)/2)=e^{i\gamma}=e^{-u}$.

If $\Pi_\rho(d)$ denotes the angle of parallelism at geometric distance $d$ in a hyperbolic plane of curvature $-1/\rho^2$, then $u=d/\rho$ and
\[
\tan\!\left(\frac{\Pi_\rho(d)}{2}\right)=e^{-d/\rho}.
\]

\section{Transformations}
\label{sec:transforms}

Throughout this section, let
\[
a>0,
\qquad b\in\mathbb R,
\qquad y:=\frac{x+b}{a}.
\]
Each transform is a change-of-variable identity on an open connected interval contained in its stated domain. Constants of integration may differ between connected components.

\subsection*{Bridge lemma for the sum radical}

\begin{lemma}[Bridge lemma]
\label{lem:bridge}
For real $y$,
\[
e^{i\arccos(iy)}=i\left(y+\sqrt{y^2+1}\right).
\]
Equivalently,
\[
-ie^{i\arccos(iy)}=y+\sqrt{y^2+1}.
\]
\end{lemma}

\begin{proof}
For real $y$, $\sqrt{1-(iy)^2}=\sqrt{1+y^2}>0$ on the principal branch, and the defining logarithmic formulas give
\[
\arcsin(iy)=i\arsinh y,
\qquad
\arsinh y=\Logop\!\left(y+\sqrt{y^2+1}\right).
\]
Therefore,
\begin{align*}
e^{i\arccos(iy)}
&=e^{i(\pi/2-\arcsin(iy))}\\
&=ie^{\arsinh y}\\
&=i\left(y+\sqrt{y^2+1}\right).
\end{align*}
\end{proof}

\subsection*{Transform 1: $\tan(\beta/2)$ and $\tan(\gamma/2)$}

Let
\[
\alpha=\arccos y,
\qquad
\beta=\arccsc y,
\qquad
\gamma=\arcsec y,
\]
and choose
\[
t=\tan\!\left(\frac{\beta}{2}\right)
=
\begin{cases}
y-\sqrt{y^2-1},&y\ge1,\\[1mm]
y+\sqrt{y^2-1},&y\le-1.
\end{cases}
\]
Then
\[
\tan\!\left(\frac{\gamma}{2}\right)=\frac{1-t}{1+t},
\qquad
x=\frac{a(t^2+1)}{2t}-b,
\qquad
dx=a\frac{t^2-1}{2t^2}\,dt.
\]
Hence
\[
\int f\!\left(x,\tan\frac{\beta}{2},\tan\frac{\gamma}{2}\right)\,dx
=
\int f\!\left(\frac{a(t^2+1)}{2t}-b,t,\frac{1-t}{1+t}\right)
 a\frac{t^2-1}{2t^2}\,dt.
\]
The antiderivative formula is applied on either open component
\[
y\in(-\infty,-1)
\qquad\text{or}\qquad
y\in(1,\infty).
\]
The algebraic formulas are well defined at $y=1$. At $y=-1$, $\tan(\gamma/2)$ is undefined; that endpoint may be used only when the complete integrand has a verified one-sided limit. The back-substitution is
\[
t=
\begin{cases}
\dfrac{x+b-\sqrt{(x+b)^2-a^2}}{a},&y\ge1,\\[3mm]
\dfrac{x+b+\sqrt{(x+b)^2-a^2}}{a},&y\le-1.
\end{cases}
\]

\subsection*{Transform 2: $\sqrt{(x+b)^2-a^2}$ and $\sqrt{(x+b-a)/(x+b+a)}$}

Use the same parameter $t$ as in Transform~1 and set
\[
\varepsilon:=\sgn(y)\in\{-1,1\}.
\]
For $|y|\ge1$,
\[
\sqrt{(x+b)^2-a^2}
=-\varepsilon a\frac{t^2-1}{2t}
=
\begin{cases}
-a\dfrac{t^2-1}{2t},&y\ge1,\\[2mm]
a\dfrac{t^2-1}{2t},&y\le-1.
\end{cases}
\]
For $y\ne-1$,
\[
\sqrt{\frac{x+b-a}{x+b+a}}=\frac{1-t}{1+t}.
\]
Consequently,
\begin{align*}
&\int f\!\left(x,\sqrt{(x+b)^2-a^2},\sqrt{\frac{x+b-a}{x+b+a}}\right)\,dx\\
&\quad=
\int f\!\left(
\frac{a(t^2+1)}{2t}-b,
-\varepsilon a\frac{t^2-1}{2t},
\frac{1-t}{1+t}
\right)
a\frac{t^2-1}{2t^2}\,dt.
\end{align*}
The antiderivative formula is applied on either open component $y\in(-\infty,-1)$ or $y\in(1,\infty)$. The radical and Jacobian formulas extend to $y=1$. The point $y=-1$ is excluded because the radical quotient is undefined there.

\subsection*{Transform 3: $\sqrt{(x+b)^2+a^2}$}

Let
\[
\theta=\arsinh y,
\qquad
s=e^\theta=y+\sqrt{y^2+1}>0.
\]
Then
\[
y=\frac{s-s^{-1}}{2},
\qquad
\sqrt{y^2+1}=\frac{s+s^{-1}}{2},
\]
so
\[
x=\frac{a(s^2-1)}{2s}-b,
\qquad
\sqrt{(x+b)^2+a^2}=\frac{a(s^2+1)}{2s},
\qquad
dx=a\frac{s^2+1}{2s^2}\,ds.
\]
Thus
\[
\int f\!\left(x,\sqrt{(x+b)^2+a^2}\right)\,dx
=
\int f\!\left(\frac{a(s^2-1)}{2s}-b,\frac{a(s^2+1)}{2s}\right)
a\frac{s^2+1}{2s^2}\,ds.
\]
This transform is valid for every real $y$. By Lemma~\ref{lem:bridge},
\[
s=-ie^{i\arccos(iy)},
\]
so the transform follows from the same principal-branch definitions through the bridge lemma.

\subsection*{Transform 4: $\tan(\psi/2)$ and $\tan(\alpha/2)$}

Let
\[
\psi=\arcsin y,
\qquad
\alpha=\arccos y,
\]
and define
\[
r:=\tan\!\left(\frac{\psi}{2}\right)
=\frac{y}{1+\sqrt{1-y^2}}.
\]
This definition also applies at $y=0$, where $r=0$. For $|y|>1$, both $\arcsin y$ and $\sqrt{1-y^2}$ are understood as the boundary values prescribed in Section~\ref{sec:branches}. For $y\ne0$, Theorem~\ref{thm:sec-sin} also gives the corresponding exponential representation in terms of $\arcsec y$; for $0<|y|<1$, that value of $\arcsec y$ is likewise a prescribed boundary value.

For $-1<y\le1$,
\[
\tan\!\left(\frac{\alpha}{2}\right)=\frac{1-r}{1+r}.
\]
Moreover,
\[
y=\frac{2r}{1+r^2},
\qquad
x=\frac{2ar}{1+r^2}-b,
\qquad
dx=2a\frac{1-r^2}{(1+r^2)^2}\,dr.
\]
Therefore,
\[
\int f\!\left(x,\tan\frac{\psi}{2},\tan\frac{\alpha}{2}\right)\,dx
=
\int f\!\left(\frac{2ar}{1+r^2}-b,r,\frac{1-r}{1+r}\right)
2a\frac{1-r^2}{(1+r^2)^2}\,dr.
\]
When $\tan(\alpha/2)$ occurs as an input, the real antiderivative formula is applied on $-1<y<1$. The algebraic formulas are well defined at $y=1$. At $y=-1$, the half-angle is undefined and only a limiting statement is available. A convenient real back-substitution on $|y|\le1$ is
\[
r=\frac{x+b}{a+\sqrt{a^2-(x+b)^2}},
\]
which gives $r=0$ at $x=-b$.

\subsection*{Transform 5: $\sqrt{a^2-(x+b)^2}$ and $\sqrt{(a+b+x)/(a-b-x)}$}

Use the same parameter $r$ as in Transform~4. For $|y|\le1$,
\[
\sqrt{a^2-(x+b)^2}=a\frac{1-r^2}{1+r^2},
\]
and, for $y<1$,
\[
\sqrt{\frac{a+b+x}{a-b-x}}=\frac{1+r}{1-r}.
\]
Hence
\begin{align*}
&\int f\!\left(x,\sqrt{a^2-(x+b)^2},\sqrt{\frac{a+b+x}{a-b-x}}\right)\,dx\\
&\quad=
\int f\!\left(
\frac{2ar}{1+r^2}-b,
a\frac{1-r^2}{1+r^2},
\frac{1+r}{1-r}
\right)
2a\frac{1-r^2}{(1+r^2)^2}\,dr.
\end{align*}
The antiderivative formula is applied on $-1<y<1$. The algebraic formulas are well defined at $y=-1$. At $y=1$, the radical quotient is undefined.

\begin{corollary}[Weierstrass substitution as a specialization of Transform~5]
\label{cor:weierstrass}
Let $J$ be an open interval on which
\[
\arcsin(\sin\omega)=\omega
\qquad\text{and}\qquad
\cos\omega>0;
\]
for example, $J\subset(-\pi/2,\pi/2)$. If $R$ is rational in its two arguments, then Transform~5 reduces
\[
I=\int R(\sin\omega,\cos\omega)\,d\omega
\]
to
\[
I=\int R\!\left(\frac{2r}{1+r^2},\frac{1-r^2}{1+r^2}\right)\frac{2}{1+r^2}\,dr,
\qquad
r=\tan\!\left(\frac{\omega}{2}\right).
\]
Thus the unit-radius specialization of Transform~5 is precisely the
corresponding local branch of the classical Weierstrass, or tangent
half-angle, substitution.\footnote{The designation \emph{Weierstrass substitution} is conventional
rather than historical; for the standard modern terminology, see
\cite{Weisstein}. The tangent half-angle substitution predates Weierstrass:
Euler used it in \emph{Institutionum calculi integralis} (1768), \S261
\cite[\S261]{Euler}, and Legendre later treated the rationalizing half-angle
parametrization in \emph{Exercices de calcul int\'egral}, Tome~2 (1817),
pp.~245--246 \cite[pp.~245--246]{Legendre}. See also \cite{Johnson} for
historical discussion of the attribution.}
\end{corollary}

\begin{proof}
Set $x=\sin\omega$. On $J$, the principal circular radical equals $\cos\omega$, so
\[
\cos\omega=\sqrt{1-x^2},
\qquad
d\omega=\frac{dx}{\sqrt{1-x^2}}.
\]
Therefore,
\[
I=\int \frac{R\!\left(x,\sqrt{1-x^2}\right)}{\sqrt{1-x^2}}\,dx.
\]
Apply Transform~5 with $a=1$ and $b=0$. It gives
\[
x=\frac{2r}{1+r^2},
\qquad
\sqrt{1-x^2}=\frac{1-r^2}{1+r^2},
\qquad
dx=\frac{2(1-r^2)}{(1+r^2)^2}\,dr.
\]
Substitution displays the cancellation explicitly:
\begin{align*}
I
&=\int
R\!\left(\frac{2r}{1+r^2},\frac{1-r^2}{1+r^2}\right)
\frac{1+r^2}{1-r^2}
\frac{2(1-r^2)}{(1+r^2)^2}\,dr\\
&=\int
R\!\left(\frac{2r}{1+r^2},\frac{1-r^2}{1+r^2}\right)
\frac{2}{1+r^2}\,dr.
\end{align*}
Because $\arcsin x=\omega$ on $J$, the Transform~5 parameter is
\[
r=\tan\!\left(\frac{1}{2}\arcsin x\right)=\tan\!\left(\frac{\omega}{2}\right).
\]
Consequently,
\[
\sin\omega=\frac{2r}{1+r^2},
\qquad
\cos\omega=\frac{1-r^2}{1+r^2},
\qquad
d\omega=\frac{2\,dr}{1+r^2},
\]
which are the classical Weierstrass formulas.
\end{proof}

\begin{proposition}[Laurent-polynomial families]
\label{prop:laurent}
Let $X=x+b$, let $P$ be a polynomial, and restrict all expressions to an open component on which the corresponding transform and principal-root representation hold and every displayed denominator is nonzero.
\begin{enumerate}[label=\arabic*.]
\item Under Transform~2, for every polynomial $P$,
\[
P(x)\sqrt{X^2-a^2}\,dx
\]
becomes a finite Laurent polynomial in $t$.
\item Under Transform~2, for every polynomial $P$,
\[
P(x)\sqrt{\frac{X-a}{X+a}}\,dx
\]
becomes a finite Laurent polynomial in $t$.
\item Under Transform~3, for every polynomial $P$, both
\[
P(x)\sqrt{X^2+a^2}\,dx,
\qquad
\frac{P(x)}{\sqrt{X^2+a^2}}\,dx
\]
become finite Laurent polynomials in $s$. Moreover, for every integer $n\ge1$, each of
\[
\frac{dx}{\left(X+\sqrt{X^2+a^2}\right)^n},
\qquad
\frac{dx}{\left(\sqrt{X^2+a^2}-X\right)^n}
\]
becomes a finite Laurent polynomial in $s$.
\item Under Transform~5, for every integer $n\ge1$,
\[
\frac{dx}{X^n\sqrt{a^2-X^2}}
\]
becomes a finite Laurent polynomial in $r$.
\item Under Transform~5, for every integer $n\ge3$,
\[
\frac{\sqrt{a^2-X^2}}{X^n}\,dx
\]
becomes a finite Laurent polynomial in $r$.
\item Under Transform~5, for every integer $n\ge2$, each of
\[
\frac{1}{X^n}\sqrt{\frac{a+X}{a-X}}\,dx,
\qquad
\frac{1}{X^n}\sqrt{\frac{a-X}{a+X}}\,dx
\]
becomes a finite Laurent polynomial in $r$.
\end{enumerate}
When $b=0$, one has $X=x$, so the last three families are precisely the corresponding families with reciprocal powers of $x$.
\end{proposition}

\begin{proof}
For Transform~2,
\[
x=\frac{a}{2}\left(t+t^{-1}\right)-b,
\]
so every polynomial $P(x)$ is a Laurent polynomial in $t$. The difference radical gives
\[
\sqrt{X^2-a^2}\,dx
=-\varepsilon\frac{a^2}{4}t^{-3}(t^2-1)^2\,dt,
\]
which is Laurent because $\varepsilon$ is constant on each component. The radical quotient cancels the factor $1+t$ in the Jacobian:
\begin{align*}
\sqrt{\frac{X-a}{X+a}}\,dx
&=\frac{1-t}{1+t}\,a\frac{t^2-1}{2t^2}\,dt\\
&=-\frac{a}{2}\left(1-t^{-1}\right)^2\,dt.
\end{align*}
Multiplication by $P(x)$ therefore gives a Laurent polynomial in both cases.

For Transform~3, write
\[
H:=\sqrt{X^2+a^2}.
\]
Then
\[
x=\frac{a}{2}\left(s-s^{-1}\right)-b,
\qquad
H=\frac{a}{2}\left(s+s^{-1}\right),
\qquad
dx=\frac{a}{2}\left(1+s^{-2}\right)\,ds.
\]
Thus every polynomial $P(x)$ is a Laurent polynomial in $s$, and
\[
H\,dx=\frac{a^2}{4}\left(s+2s^{-1}+s^{-3}\right)\,ds,
\qquad
\frac{dx}{H}=s^{-1}\,ds.
\]
Consequently, both $P(x)H\,dx$ and $P(x)H^{-1}\,dx$ are finite Laurent polynomials in $s$. In addition,
\[
X+H=as,
\qquad
H-X=as^{-1},
\]
and hence, for $n\ge1$,
\[
\frac{dx}{(X+H)^n}
=\frac{a^{1-n}}{2}\left(s^{-n}+s^{-n-2}\right)\,ds,
\]
\[
\frac{dx}{(H-X)^n}
=\frac{a^{1-n}}{2}\left(s^n+s^{n-2}\right)\,ds.
\]
Each expression is a finite Laurent polynomial.

For Transform~5,
\[
X=\frac{2ar}{1+r^2},
\qquad
\sqrt{a^2-X^2}=a\frac{1-r^2}{1+r^2},
\qquad
dx=2a\frac{1-r^2}{(1+r^2)^2}\,dr.
\]
Direct cancellation gives
\[
\frac{dx}{X^n\sqrt{a^2-X^2}}
=2^{1-n}a^{-n}r^{-n}(1+r^2)^{n-1}\,dr,
\qquad n\ge1,
\]
\[
\frac{\sqrt{a^2-X^2}}{X^n}\,dx
=2^{1-n}a^{2-n}r^{-n}(1-r^2)^2(1+r^2)^{n-3}\,dr,
\qquad n\ge3,
\]
and
\[
\frac{1}{X^n}\sqrt{\frac{a\pm X}{a\mp X}}\,dx
=2^{1-n}a^{1-n}r^{-n}(1\pm r)^2(1+r^2)^{n-2}\,dr,
\qquad n\ge2.
\]
Each right-hand side is a finite Laurent polynomial. The stated restrictions on $n$ are exactly those that prevent a residual negative power of $1+r^2$ in these formulas.
\end{proof}

\begin{remark}[Termwise integration]
If
\[
L(p)=\sum_{k=m}^{n}c_kp^k,
\]
with the convention that \(c_k=0\) for \(k\notin[m,n]\),
then, on a real component on which $p$ has fixed sign,
\[
\int L(p)\,dp
=
\sum_{\substack{k=m\\k\ne-1}}^{n}\frac{c_k}{k+1}p^{k+1}
+c_{-1}\ln|p|+C.
\]
Thus no partial-fraction decomposition is required. Examples~3, 4, and 6 reduce, respectively, to Laurent polynomials in $s$, $t$, and $r$. This is an algebraic advantage of the transforms, but it is not a universal complexity theorem for arbitrary integrands.
\end{remark}

\section{Worked examples}

\subsection*{Example 1: $\displaystyle\int\sqrt{x^2-1}\,dx$ for $x>1$}

For $x>1$, Transform~2 with $a=1$ and $b=0$ gives
\[
t=x-\sqrt{x^2-1},
\qquad
\sqrt{x^2-1}=-\frac{t^2-1}{2t},
\qquad
dx=\frac{t^2-1}{2t^2}\,dt.
\]
Thus
\begin{align*}
\sqrt{x^2-1}\,dx
&=\left(-\frac{t^2-1}{2t}\right)\left(\frac{t^2-1}{2t^2}\right)\,dt\\
&=-\frac{(t^2-1)^2}{4t^3}\,dt\\
&=-\frac{1}{4}\left(t-2t^{-1}+t^{-3}\right)\,dt.
\end{align*}
Therefore,
\begin{align*}
\int\sqrt{x^2-1}\,dx
&=-\frac{1}{4}\int\left(t-2t^{-1}+t^{-3}\right)\,dt\\
&=-\frac{t^2}{8}+\frac{1}{2}\ln t+\frac{1}{8t^2}+C.
\end{align*}
Back-substitution yields
\[
\int\sqrt{x^2-1}\,dx
=\frac{1}{2}\left(x\sqrt{x^2-1}-\ln\!\left(x+\sqrt{x^2-1}\right)\right)+C.
\]
Differentiation verifies the formula for $x>1$, and the expression extends continuously to $x=1$. On $x<-1$, a corresponding primitive is
\[
\frac{1}{2}\left(x\sqrt{x^2-1}-\ln\left|x+\sqrt{x^2-1}\right|\right)+C.
\]

\subsection*{Example 2: $\displaystyle\int\frac{dx}{x\sqrt{x^2+x}}$ for $x>0$ or $x<-1$}

On $x>0$, complete the square:
\[
x^2+x=\left(x+\frac{1}{2}\right)^2-\left(\frac{1}{2}\right)^2.
\]
Transform~2 applies with $a=b=1/2$, giving
\[
t=2x+1-2\sqrt{x^2+x},
\]
with
\[
x=\frac{(t-1)^2}{4t},
\qquad
\sqrt{x^2+x}=\frac{1-t^2}{4t},
\qquad
dx=\frac{t^2-1}{4t^2}\,dt.
\]
Hence
\[
\frac{dx}{x\sqrt{x^2+x}}
=-\frac{4}{(t-1)^2}\,dt,
\]
and therefore
\[
\int\frac{dx}{x\sqrt{x^2+x}}
=-4\int\frac{dt}{(t-1)^2}
=\frac{4}{t-1}+C
=\frac{2}{x-\sqrt{x^2+x}}+C.
\]
Direct differentiation verifies this expression on both $x>0$ and $x<-1$. On the lower component $x<-1$, the Transform~2 parameter is
\[
t=2x+1+2\sqrt{x^2+x},
\]
and the lower-component parametrization naturally gives
\[
-\frac{2}{x+\sqrt{x^2+x}}+C.
\]
On $x<-1$, the two displayed primitives differ by the constant $4$.

\subsection*{Example 3: $\displaystyle\int_0^\infty\frac{dx}{\left(x+\sqrt{1+x^2}\right)^2}$}

Set
\[
s=x+\sqrt{1+x^2}.
\]
Then
\[
x=\frac{s^2-1}{2s},
\qquad
\sqrt{1+x^2}=\frac{s^2+1}{2s},
\qquad
dx=\frac{s^2+1}{2s^2}\,ds.
\]
Moreover, $x:0\to\infty$ maps to $s:1\to\infty$. Hence
\[
\frac{dx}{\left(x+\sqrt{1+x^2}\right)^2}
=\frac{1}{2}\left(s^{-2}+s^{-4}\right)\,ds.
\]
Therefore,
\[
\int_0^\infty\frac{dx}{\left(x+\sqrt{1+x^2}\right)^2}
=\frac{1}{2}\int_1^\infty\left(s^{-2}+s^{-4}\right)\,ds
=\frac{2}{3}.
\]

\subsection*{Example 4: $\displaystyle\int\sqrt{\frac{x+1}{x+3}}\,dx$ for $x>-1$}

Matching
\[
x+1=x+b-a,
\qquad
x+3=x+b+a
\]
gives $a=1$ and $b=2$. Transform~2 yields
\[
t=x+2-\sqrt{x^2+4x+3},
\]
and
\[
\sqrt{\frac{x+1}{x+3}}=\frac{1-t}{1+t},
\qquad
dx=\frac{t^2-1}{2t^2}\,dt.
\]
Thus
\begin{align*}
\sqrt{\frac{x+1}{x+3}}\,dx
&=-\frac{(1-t)^2}{2t^2}\,dt\\
&=-\frac{1}{2}\left(1-2t^{-1}+t^{-2}\right)\,dt.
\end{align*}
Consequently,
\begin{align*}
\int\sqrt{\frac{x+1}{x+3}}\,dx
&=-\frac{1}{2}\int\left(1-2t^{-1}+t^{-2}\right)\,dt\\
&=-\frac{t}{2}+\ln t+\frac{1}{2t}+C.
\end{align*}
The transformed integrand is a Laurent polynomial and is integrated term by term, without partial fractions. Back-substitution gives
\[
\int\sqrt{\frac{x+1}{x+3}}\,dx
=\ln\!\left(x+2-\sqrt{x^2+4x+3}\right)
+\sqrt{x^2+4x+3}+C.
\]
The expression differentiates to the integrand on $x>-1$ and extends continuously to $x=-1$.

The direct rationalization
\[
u=\sqrt{\frac{x+1}{x+3}}
\]
produces
\[
\int\frac{4u^2}{(u^2-1)^2}\,du.
\]
On the stated real domain, $0\le u<1$; therefore substitutions such as $u=\coth z$ or $u=\sec\theta$ do not parametrize the relevant real range. Valid alternatives include direct partial fractions or $u=\tanh z$.

\subsection*{Example 5: $\displaystyle\int\frac{dx}{\tan\!\left(\frac12\arccsc x\right)-\tan\!\left(\frac12\arcsec x\right)}$}

On each open component
\[
x\in(-\infty,-1),
\qquad
x\in(1,\sqrt2),
\qquad\text{or}\qquad
x\in(\sqrt2,\infty),
\]
set
\[
t=\tan\!\left(\frac{1}{2}\arccsc x\right).
\]
Then
\[
\tan\!\left(\frac{1}{2}\arcsec x\right)=\frac{1-t}{1+t},
\qquad
dx=\frac{t^2-1}{2t^2}\,dt,
\]
and
\[
I=\frac{1}{2}\int\frac{(t-1)(t+1)^2}{t^2(t^2+2t-1)}\,dt.
\]
A partial-fraction decomposition is
\[
\frac{1}{2}\frac{(t-1)(t+1)^2}{t^2(t^2+2t-1)}
=-\frac{t+3}{t^2+2t-1}+\frac{3}{2t}+\frac{1}{2t^2}.
\]
Therefore,
\begin{align*}
I={}&\frac{3}{2}\ln|t|-\frac{1}{2t}
-\frac{1}{2}\ln|t^2+2t-1|\\
&-\frac{1}{\sqrt2}\ln\left|\frac{t+1-\sqrt2}{t+1+\sqrt2}\right|+C.
\end{align*}
The original denominator vanishes at
\[
t=\sqrt2-1,
\qquad x=\sqrt2,
\]
which accounts for the two positive-domain components. For back-substitution, the componentwise formulas are
\[
t=x-\sqrt{x^2-1},\qquad x>1,
\]
and
\[
t=x+\sqrt{x^2-1},\qquad x<-1.
\]
At $x=-1$, the inverse-secant half-angle is undefined. At $x=1$, the original integrand is defined, and any endpoint statement should be interpreted as a one-sided extension of a primitive from $(1,\sqrt2)$.

\subsection*{Example 6: $\displaystyle\int\frac{dx}{x^3\sqrt{4-x^2}}$ for $0<x<2$}

Transform~5 with $a=2$ and $b=0$ gives
\[
r=\frac{2-\sqrt{4-x^2}}{x},
\qquad
x=\frac{4r}{1+r^2}.
\]
The radical and Jacobian factors cancel, giving
\[
\int\frac{dx}{x^3\sqrt{4-x^2}}
=\frac{1}{32}\int\left(r^{-3}+2r^{-1}+r\right)\,dr.
\]
This is again a Laurent polynomial. Termwise integration and back-substitution give
\[
\int\frac{dx}{x^3\sqrt{4-x^2}}
=\frac{1}{16}\left[
\ln\left|\frac{2-\sqrt{4-x^2}}{x}\right|
-\frac{2\sqrt{4-x^2}}{x^2}
\right]+C.
\]
Differentiation gives the original integrand on $0<x<2$.

\subsection*{Example 7: $\displaystyle\int e^{\arccos x}\,dx$ for $-1<x<1$}

Let
\[
\alpha=\arccos x,
\qquad
t=e^{-i\alpha}=x-i\sqrt{1-x^2}.
\]
Along this path,
\[
\Logop t=-i\alpha.
\]
As $x\to-1^+$, $t\to-1$ through the lower half-plane and $\Logop t\to-i\pi$; this is a lower-bank boundary limit, not a value of the analytic logarithm at a point on its cut. Since $e^\alpha=t^i$ along the path,
\begin{align*}
\int e^{\arccos x}\,dx
&=\frac{1}{2}\int\left(t^i-t^{i-2}\right)\,dt\\
&=\frac{1}{4}t^i\left((1-i)t+(1+i)t^{-1}\right)+C.
\end{align*}
Back-substitution gives
\[
\int e^{\arccos x}\,dx
=\frac{e^{\arccos x}}{2}\left(x-\sqrt{1-x^2}\right)+C.
\]
Direct differentiation verifies the formula on $-1<x<1$, and the expression extends continuously to both endpoints.

\begin{remark}
In Example~7, for every \(\lambda\in\mathbb C\), the complex power
\(t^\lambda\) is defined using the principal logarithm:
\[
t^\lambda:=\exp\!\bigl(\lambda\Logop t\bigr).
\]
Along the path
\[
t=e^{-i\alpha},
\qquad
\alpha=\arccos x\in(0,\pi),
\]
one has
\[
\Logop t=-i\alpha.
\]
Consequently,
\[
t^i=\exp\!\bigl(i\Logop t\bigr)=e^\alpha,
\qquad
t^{\,i-2}=t^i t^{-2}.
\]
Example~7 demonstrates that the same exponential-parameter language also handles a nonalgebraic integrand involving \(\arccos x\).
\end{remark}

\section{Relation to Euler substitutions}

We compare the present transforms with the two substitutions used by
Euler in his original treatment.\footnote{Modern accounts commonly
present three ``Euler substitutions.'' However, Cie\'sli\'nski and
Jurgielewicz note that Euler himself used only two of these
substitutions in his original treatment; the now-familiar third Euler
substitution belongs to the later three-substitution classification,
whose origin they note is unclear. See
\cite[\S2.4]{CieslinskiJurgielewicz}.}

Let
\[
Q(x)=a'x^2+b'x+c',
\qquad a'\ne0,
\qquad \Delta=b'^2-4a'c',
\]
and set
\[
X=x+\frac{b'}{2a'}.
\]
Then
\[
Q(x)=a'\left(X^2-\frac{\Delta}{4a'^2}\right).
\]
For $\Delta\ne0$, define
\[
A:=\frac{\sqrt{|\Delta|}}{2|a'|}>0.
\]
Restricting first to real-valued principal radicals, the nondegenerate cases are
\[
\begin{array}{lll}
\text{difference form:}&\Delta>0,\ a'>0,
&\sqrt{Q(x)}=\sqrt{a'}\sqrt{X^2-A^2},\quad |X|\ge A,\\[1mm]
\text{circular form:}&\Delta>0,\ a'<0,
&\sqrt{Q(x)}=\sqrt{|a'|}\sqrt{A^2-X^2},\quad |X|\le A,\\[1mm]
\text{sum form:}&\Delta<0,\ a'>0,
&\sqrt{Q(x)}=\sqrt{a'}\sqrt{X^2+A^2}.
\end{array}
\]
If $\Delta<0$ and $a'<0$, then $Q(x)<0$ for all real $x$, and its prescribed principal boundary value is
\[
\sqrt{Q(x)}=i\sqrt{|a'|}\sqrt{X^2+A^2}.
\]
If $\Delta=0$, then $Q(x)=a'X^2$ and
\[
\sqrt{Q(x)}=
\begin{cases}
\sqrt{a'}\,|X|,&a'>0,\\[1mm]
i\sqrt{|a'|}\,|X|,&a'<0.
\end{cases}
\]

\subsection*{Difference form and Euler's first substitution}

Let
\[
S=\sqrt{X^2-A^2},
\qquad
t=
\begin{cases}
\dfrac{X-S}{A},&X\ge A,\\[2mm]
\dfrac{X+S}{A},&X\le-A.
\end{cases}
\]
Then
\[
X=\frac{A}{2}\left(t+t^{-1}\right),
\qquad
S=-\sgn(X)\frac{A}{2}\left(t-t^{-1}\right).
\]
The standard Euler combination
\[
u_E:=\sqrt{Q(x)}-\sqrt{a'}\,X
=\sqrt{a'}(S-X)
\]
satisfies
\[
u_E=
\begin{cases}
-A\sqrt{a'}\,t,&X\ge A,\\[1mm]
-A\sqrt{a'}\,t^{-1},&X\le-A.
\end{cases}
\]
Thus Transform~2 recovers Euler's first rationalization after an explicit scaling and a component-dependent reciprocal reparametrization.

\subsection*{Circular form and Euler's second substitution}

Let
\[
K=\sqrt{A^2-X^2},
\qquad
r=\frac{X}{A+K}.
\]
Then
\[
X=\frac{2Ar}{1+r^2},
\qquad
K=A\frac{1-r^2}{1+r^2}.
\]
If the normalized Euler parameter is chosen as
\[
t_E:=\frac{K-A}{X},
\]
with its value at $X=0$ understood by continuity, then
\[
t_E=-r.
\]
Accordingly, Transform~5 agrees with this sign convention for Euler's second substitution. Other textbook conventions differ by a sign.

\subsection*{Sum form}

Transform~3 gives
\[
X=A\frac{s^2-1}{2s},
\qquad
\sqrt{X^2+A^2}=A\frac{s^2+1}{2s}.
\]
The reciprocal parameter $1/s$, with an optional scale factor, yields the standard Euler-type rational forms for the sum radical.

\section{Conclusion}

The Unified Substitution Method places several classical rationalizations in a common branch-explicit framework. Its three parameters generate five transformations in which the variable, radical or half-angle expression, and differential become rational functions of a single parameter. The branch information is essential: the sum-radical parameter is globally positive, the circular parameter is continuous under the adopted boundary-value convention, and the difference-radical parameter $t$, defined through the inverse-cosecant half-angle, has component-dependent forms and no continuous extension through zero. Algebraic reconstruction alone therefore does not determine radical signs, logarithmic branches, or constants of integration.

For the families in Proposition~\ref{prop:laurent}, the transformed integrands are finite Laurent polynomials and can be integrated term by term. The method also contains the corresponding local branch of the Weierstrass substitution, recovers Euler rationalizations after the necessary sign and scaling choices, and yields Lobachevsky's formula
\[
\tan\!\left(\frac{\Pi(u)}{2}\right)=e^{-u}.
\]
Under the author-reported protocol, all tested instances for Transforms~2--5 passed the reported derivative/numerical checks, as did 94 of 100 instances for Transform~1. The reported experiments indicate reduced expression swell or improved completion coverage for selected targeted families, but they do not establish superiority for arbitrary integrals. The principal contribution is structural: rationalization and the associated branch conditions are specified within one framework.

\appendix

\section{Complementary half-angle identities and parameter continuation}
\label{app:continuation}

Throughout this appendix,
\[
y=\frac{x+b}{a},
\qquad a>0,
\]
and the branch conventions of Section~\ref{sec:branches} remain in force.

\begin{theorem}[Interior continuation of the $t$-half-angle]
\label{thm:interior-t}
Let
\[
\alpha=\arccos y,
\qquad
\beta=\arccsc y,
\qquad
0<|y|\le1.
\]
For $0<|y|<1$, $\beta$ denotes the boundary value prescribed in Section~\ref{sec:branches}; at $y=\pm1$, it denotes the continuous endpoint value. Then
\[
\tan\!\left(\frac{\beta}{2}\right)=
\begin{cases}
e^{i\alpha}=y+i\sqrt{1-y^2},&-1\le y<0,\\[1mm]
e^{-i\alpha}=y-i\sqrt{1-y^2},&0<y\le1.
\end{cases}
\]
Moreover,
\[
\lim_{y\to0^-}\tan\!\left(\frac{1}{2}\arccsc y\right)=i,
\qquad
\lim_{y\to0^+}\tan\!\left(\frac{1}{2}\arccsc y\right)=-i.
\]
\end{theorem}

\begin{proof}
Since $\beta=\arcsin(1/y)$,
\[
\sin\beta=\frac{1}{y},
\qquad
\cos\beta=i\frac{\sqrt{1-y^2}}{|y|}.
\]
Using
\[
\tan\!\left(\frac{\beta}{2}\right)=\frac{1-\cos\beta}{\sin\beta}
\]
gives the two component formulas. The one-sided limits follow immediately.
\end{proof}

\begin{theorem}[Four-component form of the $t$-parameter]
\label{thm:four-t}
For $y\ne0$, define
\[
t:=\tan\!\left(\frac{1}{2}\arccsc y\right).
\]
On $0<|y|<1$, $\arccsc y$ denotes the boundary value prescribed in Section~\ref{sec:branches}. Then
\[
t=
\begin{cases}
y+\sqrt{y^2-1},&y\le-1,\\[1mm]
y+i\sqrt{1-y^2},&-1<y<0,\\[1mm]
y-i\sqrt{1-y^2},&0<y<1,\\[1mm]
y-\sqrt{y^2-1},&y\ge1.
\end{cases}
\]
On every component,
\[
y=\frac{1}{2}\left(t+t^{-1}\right),
\qquad
x=\frac{a}{2}\left(t+t^{-1}\right)-b,
\qquad
dx=a\frac{t^2-1}{2t^2}\,dt.
\]
The parameter agrees continuously at $y=\pm1$, but it has no continuous extension through $y=0$ that preserves the stated half-angle definition.
\end{theorem}

\begin{proof}
The exterior formulas follow from Theorem~\ref{thm:cos-csc}; the interior formulas follow from Theorem~\ref{thm:interior-t}. In every case, $t=e^{\pm i\alpha}$, so
\[
y=\cos\alpha=\frac{1}{2}\left(t+t^{-1}\right).
\]
Differentiation gives the Jacobian. The one-sided limits at zero are $i$ and $-i$.
\end{proof}

\begin{theorem}[Global form of the $r$-parameter]
\label{thm:global-r}
For every real $y$, define
\[
r:=\tan\!\left(\frac{1}{2}\arcsin y\right)
=\frac{y}{1+\sqrt{1-y^2}},
\qquad r(0)=0.
\]
For $|y|>1$, both $\arcsin y$ and $\sqrt{1-y^2}$ denote the boundary values prescribed in Section~\ref{sec:branches}. Then $r:\mathbb R\to\mathbb C$ is continuous, and
\[
y=\frac{2r}{1+r^2},
\qquad
x=\frac{2ar}{1+r^2}-b,
\qquad
dx=2a\frac{1-r^2}{(1+r^2)^2}\,dr.
\]
For $y\ne0$, the exponential representations are
\[
r=
\begin{cases}
e^{i\arcsec y},&y<-1,\\[1mm]
e^{-i\arcsec y},&-1\le y<0,\\[1mm]
e^{i\arcsec y},&0<y\le1,\\[1mm]
e^{-i\arcsec y},&y>1.
\end{cases}
\]
Here $\arcsec y$ is a real principal value for $|y|\ge1$ and the prescribed boundary value for $0<|y|<1$.
\end{theorem}

\begin{proof}
The global algebraic formula follows from
\[
\tan\!\left(\frac{\psi}{2}\right)=\frac{\sin\psi}{1+\cos\psi},
\qquad \psi=\arcsin y.
\]
Its value tends to zero from both sides at $y=0$, to $1$ as $y\to1^\pm$, and to $-1$ as $y\to-1^\pm$. Hence the displayed algebraic formula defines a continuous map from $\mathbb R$ to $\mathbb C$. Solving the quadratic relation gives $y=2r/(1+r^2)$; for finite real $y$, the denominator $1+r^2$ is nonzero. Differentiation gives the Jacobian. The exponential signs are those in Theorem~\ref{thm:sec-sin}.
\end{proof}

\subsection*{Principal-root sign tables on complementary components}

Continuation of a parameter does not imply that every principal square root in a transform retains the same rational representative.

For the $t$-parameter of Theorem~\ref{thm:four-t},
\[
\sqrt{y^2-1}=-\sgn(y)\frac{t^2-1}{2t},
\qquad y\ne0,
\]
where the square root is interpreted according to Section~\ref{sec:branches}. The inverse-secant half-angle satisfies
\[
\tan\!\left(\frac{1}{2}\arcsec y\right)=\frac{1-t}{1+t}
\]
wherever both sides are defined. However, the principal radical quotient from Transform~2 satisfies
\[
\sqrt{\frac{y-1}{y+1}}=
\begin{cases}
-\dfrac{1-t}{1+t},&-1<y<0,\\[2mm]
\dfrac{1-t}{1+t},&y\in(-\infty,-1)\cup(0,\infty),
\end{cases}
\]
with $y=-1$ excluded.

For the $r$-parameter of Theorem~\ref{thm:global-r},
\[
\sqrt{1-y^2}=\frac{1-r^2}{1+r^2}
\]
for every real $y$, using the prescribed boundary value when $|y|>1$. Also,
\[
\tan\!\left(\frac{1}{2}\arccos y\right)=\frac{1-r}{1+r}
\]
wherever both sides are defined. The principal radical quotient from Transform~5 satisfies
\[
\sqrt{\frac{1+y}{1-y}}=
\begin{cases}
\dfrac{1+r}{1-r},&y<1,\\[2mm]
-\dfrac{1+r}{1-r},&y>1,
\end{cases}
\]
and is undefined at $y=1$.

Continuation of a substitution parameter does not by itself imply continuation of the resulting antiderivative across a singularity or branch point. Antiderivative formulas are therefore understood componentwise, with the relevant logarithmic and radical branches fixed on each connected interval; different components may require different constants of integration.

\subsection*{Why Transform~3 requires no real component split}

For real $y$,
\[
s=y+\sqrt{y^2+1}
\]
is positive and defines a smooth bijection
\[
\mathbb R\longrightarrow(0,\infty),
\qquad
y\longmapsto s,
\]
with inverse $y=(s-s^{-1})/2$. Thus the same positive parameter and the same positive radical $\sqrt{y^2+1}$ apply on the entire real axis. This statement concerns only that radical and parameter; other functions in an integrand may introduce independent singularities or cuts.

\section{Binomial sum-and-difference identities}
\label{app:binomial}

Paired positive and negative powers of the substitution parameters can often be combined before back-substitution. The relevant binomial identities have different forms for the difference-radical, sum-radical, and circular parameters.

\subsection*{The difference-radical parameter $t$}

For $|y|\ge1$, let
\[
t=
\begin{cases}
y-\sqrt{y^2-1},&y\ge1,\\[1mm]
y+\sqrt{y^2-1},&y\le-1.
\end{cases}
\]
Set
\[
v=y,
\qquad
w=\sqrt{y^2-1}.
\]
Since
\[
(v-w)(v+w)=1,
\]
the unordered pair $\{v-w,v+w\}$ equals $\{t,t^{-1}\}$.

For an integer $n\ge1$, define
\[
D_n^{(t)}(y):=(v-w)^n-(v+w)^n.
\]
The binomial theorem gives
\[
D_n^{(t)}(y)
=-2\sum_{j=0}^{\lfloor(n-1)/2\rfloor}
\binom{n}{2j+1}
y^{n-2j-1}\left(\sqrt{y^2-1}\right)^{2j+1}.
\]
Because the chosen small root changes from $v-w$ to $v+w$ between the two exterior components,
\[
t^n-t^{-n}=\sigma D_n^{(t)}(y),
\qquad
\sigma=
\begin{cases}
1,&y\ge1,\\
-1,&y\le-1.
\end{cases}
\]
For example,
\[
D_1^{(t)}=-2\sqrt{y^2-1},
\qquad
D_2^{(t)}=-4y\sqrt{y^2-1},
\]
\[
D_3^{(t)}=-2(4y^2-1)\sqrt{y^2-1}.
\]
The factor $\sigma$ must not be omitted during back-substitution.

For an integer $n\ge0$, define the companion sum
\[
S_n^{(t)}(y):=(v-w)^n+(v+w)^n.
\]
Then
\[
S_n^{(t)}(y)
=2\sum_{j=0}^{\lfloor n/2\rfloor}
\binom{n}{2j}
y^{n-2j}\left(\sqrt{y^2-1}\right)^{2j}.
\]
Addition is unaffected by interchanging $v-w$ and $v+w$, so
\[
t^n+t^{-n}=S_n^{(t)}(y)
\]
on both exterior components, with no factor $\sigma$. Thus paired powers $t^n\pm t^{-n}$ can be converted directly to expressions in $y$ and $\sqrt{y^2-1}$. The difference contains a factor $\sqrt{y^2-1}$, whereas the sum becomes a polynomial in $y$ after replacing
\[
\left(\sqrt{y^2-1}\right)^{2j}=(y^2-1)^j.
\]

\subsection*{The sum-radical parameter $s$}

For every real $y$, let
\[
H:=\sqrt{y^2+1},
\qquad
s=y+H>0.
\]
Because
\[
(H+y)(H-y)=1,
\]
one has
\[
s=H+y,
\qquad
s^{-1}=H-y.
\]
For an integer $n\ge1$, define
\[
D_n^{(s)}(y):=s^n-s^{-n}.
\]
Applying the binomial theorem to $(H+y)^n-(H-y)^n$ gives
\[
D_n^{(s)}(y)
=2\sum_{j=0}^{\lfloor(n-1)/2\rfloor}
\binom{n}{2j+1}y^{2j+1}H^{n-2j-1}.
\]
Equivalently,
\[
s^n-s^{-n}
=2\sum_{j=0}^{\lfloor(n-1)/2\rfloor}
\binom{n}{2j+1}y^{2j+1}\left(\sqrt{y^2+1}\right)^{n-2j-1}.
\]
No component-dependent sign occurs because the same positive parameter $s=y+\sqrt{y^2+1}$ is used on the entire real axis. For example,
\[
s-s^{-1}=2y,
\qquad
s^2-s^{-2}=4y\sqrt{y^2+1},
\]
\[
s^3-s^{-3}=2y(4y^2+3).
\]
For an integer $n\ge0$, define
\[
S_n^{(s)}(y):=s^n+s^{-n}.
\]
Then
\[
S_n^{(s)}(y)
=2\sum_{j=0}^{\lfloor n/2\rfloor}
\binom{n}{2j}y^{2j}H^{n-2j},
\]
or explicitly,
\[
s^n+s^{-n}
=2\sum_{j=0}^{\lfloor n/2\rfloor}
\binom{n}{2j}y^{2j}\left(\sqrt{y^2+1}\right)^{n-2j}.
\]
For example,
\[
s+s^{-1}=2\sqrt{y^2+1},
\qquad
s^2+s^{-2}=4y^2+2.
\]
Consequently, every paired expression $s^n\pm s^{-n}$ reduces to a polynomial in $y$ and $\sqrt{y^2+1}$. More precisely, if $n$ is odd, the sum contains a factor $\sqrt{y^2+1}$ and the difference is a polynomial in $y$; if $n$ is even, the roles are reversed.

\subsection*{The circular parameter $r$}

For every real $y$, let
\[
c:=\sqrt{1-y^2},
\qquad
r=\frac{y}{1+c},
\]
where the square root is interpreted according to Section~\ref{sec:branches}. Thus $c$ is the nonnegative real radical for $|y|\le1$ and the prescribed boundary value for $|y|>1$.

For $y\ne0$, rationalization gives
\[
r=\frac{1-c}{y},
\qquad
r^{-1}=\frac{1+c}{y}.
\]
Equivalently, with
\[
v=\frac{1}{y},
\qquad
w=\frac{c}{y},
\]
one has
\[
r=v-w,
\qquad
r^{-1}=v+w,
\qquad
(v-w)(v+w)=1.
\]
For an integer $n\ge1$, define
\[
D_n^{(r)}(y):=r^n-r^{-n}.
\]
The binomial theorem gives
\[
D_n^{(r)}(y)
=-2\sum_{j=0}^{\lfloor(n-1)/2\rfloor}
\binom{n}{2j+1}v^{n-2j-1}w^{2j+1}.
\]
Substituting $v=1/y$ and $w=c/y$ yields
\[
r^n-r^{-n}
=-\frac{2}{y^n}
\sum_{j=0}^{\lfloor(n-1)/2\rfloor}
\binom{n}{2j+1}
\left(\sqrt{1-y^2}\right)^{2j+1},
\qquad y\ne0.
\]
For example,
\[
r-r^{-1}=-\frac{2\sqrt{1-y^2}}{y},
\qquad
r^2-r^{-2}=-\frac{4\sqrt{1-y^2}}{y^2},
\]
\[
r^3-r^{-3}=-\frac{2(4-y^2)\sqrt{1-y^2}}{y^3}.
\]
For an integer $n\ge0$, define
\[
S_n^{(r)}(y):=r^n+r^{-n}.
\]
For $n\ge1$, the corresponding binomial sum is
\[
r^n+r^{-n}
=\frac{2}{y^n}
\sum_{j=0}^{\lfloor n/2\rfloor}
\binom{n}{2j}
\left(\sqrt{1-y^2}\right)^{2j},
\qquad y\ne0.
\]
In particular,
\[
r+r^{-1}=\frac{2}{y},
\qquad
r^2+r^{-2}=\frac{4}{y^2}-2,
\]
\[
r^3+r^{-3}=\frac{2(4-3y^2)}{y^3}.
\]
Because
\[
\left(\sqrt{1-y^2}\right)^{2j}=(1-y^2)^j,
\]
the sum $r^n+r^{-n}$ is a rational function of $y$, whereas the difference $r^n-r^{-n}$ is $\sqrt{1-y^2}$ times a rational function of $y$. These identities hold for every real $y\ne0$ under the branch and boundary-value conventions of Section~\ref{sec:branches}. At $y=0$, one has $r=0$, so expressions containing negative powers of $r$ are undefined.

The three pairs of formulas permit direct reconstruction of paired Laurent powers:
\[
t^n\pm t^{-n},
\qquad
s^n\pm s^{-n},
\qquad
r^n\pm r^{-n}.
\]
The $t$-difference formula requires the component factor $\sigma$; the $s$-formulas are global on the real axis; and the $r$-formulas apply for $y\ne0$, with the adopted continuation of $\sqrt{1-y^2}$ outside $[-1,1]$.

\subsection*{Example 8: Reconstruction using \(n=1,2,3\)}

Suppose that an antiderivative obtained in \(t\)-space has the form
\[
I=
-\frac{1}{64}
\left\{
\frac{1}{3}\left(t^3-t^{-3}\right)
+\left(t^2-t^{-2}\right)
-\left(t-t^{-1}\right)
-4\ln|t|
\right\}
+C.
\]
From the difference formulas above,
\[
D^{(t)}_1(y)=-2\sqrt{y^2-1},
\qquad
D^{(t)}_2(y)=-4y\sqrt{y^2-1},
\]
and
\[
D^{(t)}_3(y)
=-2(4y^2-1)\sqrt{y^2-1}.
\]
Since
\[
t^n-t^{-n}=\sigma D^{(t)}_n(y),
\qquad
\sigma=
\begin{cases}
1, & y\ge 1,\\
-1, & y\le -1,
\end{cases}
\]
substitution gives
\[
I=
-\frac{1}{64}
\left\{
\sigma\sqrt{y^2-1}
\left[
-\frac{2}{3}(4y^2-1)-4y+2
\right]
-4\ln|t|
\right\}
+C.
\]
Equivalently,
\[
I=
\frac{\sigma}{48}
(2y^2+3y-2)\sqrt{y^2-1}
+\frac{1}{16}\ln|t|
+C.
\]

Finally, return to the original variable by using
\[
y=\frac{x+b}{a},
\]
and
\[
t=
\begin{cases}
y-\sqrt{y^2-1}
=
\dfrac{x+b-\sqrt{(x+b)^2-a^2}}{a},
& y\ge 1,\\[10pt]
y+\sqrt{y^2-1}
=
\dfrac{x+b+\sqrt{(x+b)^2-a^2}}{a},
& y\le -1.
\end{cases}
\]
Thus the component-dependent factor \(\sigma\) must be retained during
back-substitution. The square roots and logarithm are interpreted on the
corresponding real component according to the branch and domain conventions
of Sections~2 and~5.

\section{Computational benchmark and reproducibility}
\label{app:benchmark}

The counts and timings in this appendix are reproduced from the public post and linked report. No independent replication of the benchmark runs is claimed here.

The public benchmark post \cite{GarciaBenchmark} links a companion report covering five targeted benchmark families and the displayed Wolfram Language material used for the USM and native \texttt{Integrate} workflows. The linked report defines ``cold time'' as the first timed call after \texttt{ClearSystemCache[]} and measures ``warm time'' using \texttt{RepeatedTiming}; it labels an output ``monster'' when its \texttt{ByteCount} is at least $10{,}000$. These timings are local comparative measurements, not machine-independent constants.

\begin{table}[htbp]
\centering
\caption{High-level status reported in the public benchmark post and linked report.}
\label{tab:benchmark}
\begin{tabular}{lrrrrr}
\toprule
Family & USM run & USM non-timeout & USM verified & Native run & Native closed\\
\midrule
Transform 1 & 100 & 95 & 94 & 10 & 6\\
Transform 2 & 100 & 100 & 100 & 10 & 4\\
Transform 3 & 20 & 20 & 20 & 20 & 20\\
Transform 4 & 20 & 20 & 20 & 10 & 6\\
Transform 5 & 20 & 20 & 20 & 10 & 5\\
\bottomrule
\end{tabular}
\end{table}

Table~\ref{tab:benchmark} reproduces only the high-level completion counts. A reported USM verification means that the returned expression was differentiated and compared numerically with the original integrand at the benchmark's selected test points. Such finite-point agreement is a diagnostic, not a proof of a branch identity.

The report presents a differentiated picture rather than a single global performance claim. On the mutually closed Transform~1 and Transform~4 subsets, the USM expressions were reported to be substantially smaller, and the reported cold and warm times were lower. Transform~2 and Transform~5 mainly showed greater completion coverage under the stated timeout convention. For Transform~3, both methods closed all examples; USM had the lower reported cold time, whereas native \texttt{Integrate} had the lower warm time. Complete aggregate tables and the displayed code listings appear in the report linked from the public post \cite{GarciaBenchmark}.

The benchmark has several limitations:
\begin{enumerate}[label=\arabic*.]
\item The families were designed to match the corresponding transforms and are not a neutral sample of symbolic integrals.
\item Sample sizes differ among transforms and between the USM and native runs.
\item The report did not apply the same derivative-verification table to the native outputs.
\item The linked report states aggregate results larger than the example
counts in its displayed Transform~1 and Transform~2 listings; the displayed
listings therefore do not by themselves reproduce the reported aggregate counts.
\item The Transform~1 report records one non-timeout verification exception caused by an unnormalized negative \texttt{ArcCsc} argument and supplies the required sign correction.
\end{enumerate}
Accordingly, under the reported local protocol, the explicit rationalizations
gave high completion rates and controlled expression size on the targeted
families. The benchmark does not establish that USM is faster than \texttt{Integrate} for arbitrary integrals.


\begin{thebibliography}{99}

\bibitem{AbramowitzStegun}
Milton Abramowitz and Irene A. Stegun, eds., \emph{Handbook of Mathematical Functions with Formulas, Graphs, and Mathematical Tables}, National Bureau of Standards, Washington, DC, 1964. Available at \url{https://archive.org/details/handbookofmathem1964abra} (accessed October 27, 2025).

\bibitem{Ahlfors}
Lars V. Ahlfors, \emph{Complex Analysis: An Introduction to the Theory of Analytic Functions of One Complex Variable}, 3rd ed., McGraw--Hill, New York, 1979.

\bibitem{BrownChurchill}
James Ward Brown and Ruel V. Churchill, \emph{Complex Variables and Applications}, 8th ed., McGraw--Hill, New York, 2009.

\bibitem{Coxeter}
H. S. M. Coxeter, \emph{Introduction to Geometry}, 2nd ed., John Wiley \& Sons, New York, 1969, \S16.3, pp.~291--295.

\bibitem{CieslinskiJurgielewicz}
Jan L. Cie\'sli\'nski and Maciej Jurgielewicz,
``On geometric interpretation of Euler's substitutions,''
\emph{arXiv preprint} arXiv:2310.12160 [math.HO], 2023.
Available at \url{https://arxiv.org/abs/2310.12160}
(accessed August 9, 2026).

\bibitem{DLMF}
NIST Digital Library of Mathematical Functions, ``Inverse Trigonometric Functions,'' \S4.23. Available at \url{https://dlmf.nist.gov/4.23} (accessed July 19, 2026).

\bibitem{Euler}
Leonhard Euler,
\emph{Institutionum calculi integralis}, Volumen primum,
Impensis Academiae Imperialis Scientiarum, Petropoli, 1768.
Euler Archive, E342; the PDF presently supplied by the archive is
the third edition (1824).
Available at
\url{https://scholarlycommons.pacific.edu/euler-works/342/}
(accessed January 1, 2026).

\bibitem{GarciaMollweide}
Emmanuel Antonio Jos\'e Garc\'ia, ``A Generalization of Mollweide's Formula (Rather Newton's),'' \emph{Matinf} 5 (9--10) (2022), 19--22. Available at \url{http://matinf.upit.ro/MATINF9_10/files/downloads/Revista_2022.pdf#page=19} (accessed October 27, 2025).

\bibitem{GarciaRoots}
Emmanuel Jos\'e Garc\'ia, ``A Family of Trigonometric Formulas for the Roots of Quadratic Equations,'' \emph{GeoDom} (blog), February 28, 2024. Available at \url{https://geometriadominicana.blogspot.com/2024/02/a-family-of-trigonometric-formulas-for.html} (accessed October 27, 2025).

\bibitem{GarciaBenchmark}
Emmanuel Jos\'e Garc\'ia, ``Benchmarking the USM Against Mathematica Integrate, Part 3,'' \emph{GeoDom} (blog), July 8, 2026. Available at \url{https://geometriadominicana.blogspot.com/2026/07/benchmarking-usm-against-mathematica.html} (accessed July 19, 2026).

\bibitem{Hardy}
G. H. Hardy, \emph{The Integration of Functions of a Single Variable}, Cambridge University Press, Cambridge, 1905; reissued 1916. Available at \url{https://archive.org/details/integrationoffun00hardrich} (accessed October 27, 2025).

\bibitem{Johnson}
Warren P. Johnson, ``Down with Weierstrass!,'' \emph{The American Mathematical Monthly} 127 (7) (2020), 649--653. \url{https://doi.org/10.1080/00029890.2020.1763122} (accessed January 1, 2026).

\bibitem{Kishan}
Hari Kishan, \emph{Integral Calculus}, Atlantic Publishers \& Distributors, New Delhi, 2005.

\bibitem{Legendre}
Adrien-Marie Legendre,
\emph{Exercices de calcul int\'egral sur divers ordres de
transcendantes et sur les quadratures}, Tome~2,
Courcier, Paris, 1817, pp.~245--246.
Available at
\url{https://archive.org/details/exercicescalculi02legerich}
(accessed January 1, 2026).

\bibitem{Lobachevsky}
Nikolai I. Lobachevsky, \emph{Geometrical Researches on the Theory of Parallels}, translated by George Bruce Halsted, Open Court Publishing Company, Chicago, 1914, especially \S36. Available at \url{https://archive.org/details/geometricalresea00lobauoft} (accessed July 25, 2026).

\bibitem{Piskunov}
N. Piskunov, \emph{Differential and Integral Calculus}, Mir Publishers, Moscow, 1969. Available at \url{https://archive.org/details/n.-piskunov-differential-and-integral-calculus-mir-1969} (accessed October 27, 2025).

\bibitem{Weisstein}
Eric W. Weisstein, ``Weierstrass Substitution,'' \emph{MathWorld---A Wolfram Web Resource}. Available at \url{https://mathworld.wolfram.com/WeierstrassSubstitution.html} (accessed January 1, 2026).

\end{thebibliography}
\end{document}